\documentstyle[fleqn,12pt]{article}
\begin{document}
\pagenumbering{arabic}
\setcounter{page}{1}
\pagestyle{plain}
\baselineskip=18pt

\thispagestyle{empty}
\rightline{MSUMB 98-4, November 1998} 
\vspace{1.4cm}

\begin{center}
{\Large\bf Differential Geometry of $GL_{p,q}(1\vert 1)$ } 
\end{center}

\vspace{1cm}
\begin{center} Salih Celik 
\footnote{{\bf E-mail}: scelik@fened.msu.edu.tr}\\ 
Mimar Sinan University, Department of Mathematics, \\
80690 Besiktas, Istanbul, TURKEY.
\end{center}

\vspace{1.5cm}
{\bf Abstract}

We introduce a construction of the differential calculus 
on the quantum supergroup GL$_{p,q}(1\vert 1)$. We obtain two 
differential calculi respectively associated with the left and 
the right Cartan-Maurer one-forms. We also obtain the quantum 
Lie algebra of GL$_{p,q}(1\vert 1)$. Although all of the structures 
we obtain are derived without an R-matrix they nevertheless can be 
expressed using an $R$-matrix. 

\vfill\eject
\noindent
{\bf I. INTRODUCTION}

In the last few years, the theory of quantum (super) groups$^1$ has 
occured as a natural generalization of the notion of Lie groups. In 
other words, quantum (super) groups are particular deformations of 
Lie (super) groups. They are algebraic structures depending on one 
(or more) continuous parameter $q$. We have a stardard Lie (super) 
group for a particular value of the deformation parameter. Quantum 
(super) groups present the examples of (graded) Hopf algebras. They 
have found application in diverse areas of physics and mathematics$^2$. 

Quantum (super) groups can be realized on a quantum (super) space in which 
coordinates are noncommuting$^3$. Recently the differential calculus on 
noncommutative (super) space has been intensively studied both by 
mathematicians and mathematical physicists. There is much activity 
in differential geometry on quantum groups. 
Throughout the recent development of differential calculus on the 
quantum groups two principal concepts are readily seen. First of them, 
formulated by Woronowicz$^4$, is known as bicovariant differential calculus 
on the quantum groups. Another concept, introduced by Woronowicz$^5$ and 
Schirmacher {\it et al} $^6$ proceeds from the requirement of a calculus only. 
There are many papers in this field$^7$ (and references therein). 
We shall consider the second concept. 

The differential calculus on the quantum supergroups involves functions 
on the supergroup, differentials and differential forms. In ref. 8 
a right-invariant differential calculus on the quantum supergroup 
GL$_q(1\vert 1)$ has been constructed in a different way which may be 
considered as an alternative to the approach proposed earlier by 
Schmidke {\it et al} $^9$. It is necessary to point out that in the work of 
ref. 8, the generating elements of Gr$_q(1\vert 1)$ $^{10}$ have been 
interpreted as differentials of coordinate functions on the quantum supergroup 
GL$_q(1\vert 1)$ $(\hat T = \delta_R T$, in the ref. 8 notation). Therefore, 
the relations of Gr$_q(1\vert 1)$ have been used in the beginning. 
In this work, we shall construct a left-(right-) invariant 
differential calculus on the two parametric quantum supergroup, 
GL$_{p,q}(1\vert 1)$. All the relations will be obtained via some direct 
calculations and it will be shown that all relations can be rederived using 
the R-matrix approach. The differential structure obtained turns out to be a 
differential (graded) algebra. 

Let us briefly discuss the content of the paper. In the second 
section the basic notations of the Hopf algebra structure on the 
quantum supergroup $GL_{p,q}(1\vert 1)$ are introduced. In the 
third section we shall obtain the commutation relations for the 
group parameters (the matrix elements) and their differentials 
so we have a left differential algebra. We show that the 
obtained left (right) differential algebra (extended algebra) 
has a Hopf algebra structure. We also describe the quantum Lie 
algebra for the left vector fields (Lie superalgebra generators) 
for GL$_{p,q}(1\vert 1)$ and derive the commutation relations between 
the group parameters and the algebra generators. In the following 
section we propose a right differential calculus on the 
GL$_{p,q}(1\vert 1)$ and we again obtain the quantum Lie 
superalgebra of GL$_{p,q}(1\vert 1)$. In the next section we show that the 
found results can be written with help of a matrix $\hat R$. The classical 
limit $p, q\longrightarrow 1$ of the left (right) differential calculus 
gives the undeformed differential calculus. 

\noindent
{\bf II. REVIEW OF GL$_{p,q}(1\vert 1)$ }

Elementary properties of quantum supergroup $GL_{p,q}(1\vert 1)$ are 
described in$^{11}$. We state briefly the properties we are going to 
need in this work. 

The quantum supergroup $GL_{p,q}(1\vert 1)$ is defined by the matrices 
of the form 
$$T = \left(\matrix{ a & \beta \cr \gamma & d \cr} \right) = (T^i_j)\eqno(1)$$ 
where the matrix elements satisfy the following commutation 
relations$^{3,11}$ 
$$ a \beta = q \beta a, \qquad d \beta = q \beta d, $$
$$ a \gamma = p \gamma a, \qquad d \gamma = p \gamma d, \eqno(2)$$
$$ \beta \gamma + p q^{-1} \gamma \beta = 0, \qquad \beta^2 = 0 = \gamma^2, $$
$$ a d = d a + (p - q^{-1}) \gamma \beta. $$
Let us denote the algebra generated by the elements $a$, $\beta$, $\gamma$, 
$d$ with the relations (2) by ${\cal A}$. We know that the algebra 
${\cal A}$ is a (graded) Hopf algebra with the following co-structures: 

{\bf 1.} The usual coproduct 
$$\Delta : {\cal A} \longrightarrow {\cal A} \otimes {\cal A}, \qquad 
   \Delta(T^i_j) = T^i_k \otimes T^k_j,  \eqno(3)$$

{\bf 2.} the counit 
$$ \varepsilon : {\cal A} \longrightarrow {\cal C}, \qquad 
    \varepsilon(T^i_j) = \delta^i_j, \eqno(4)$$

{\bf 3.} the coinverse (antipode) $S : {\cal A} \longrightarrow {\cal A}$
$$S(T) = 
  \left(\matrix{ 
    a^{-1} + a^{-1} \beta d^{-1} \gamma a^{-1} & - a^{-1} \beta d^{-1} \cr 
    - d^{-1} \gamma a^{-1} & d^{-1} + d^{-1} \gamma a^{-1} \beta d^{-1} \cr }
     \right).  \eqno(5)$$
It is not difficult to verify the following properties of the 
co-structures: 
$$(\Delta \otimes id) \circ \Delta = (id \otimes \Delta) \circ \Delta, 
\eqno(6\mbox{a})$$
$$\mu \circ (\varepsilon \otimes id) \circ \Delta 
  = \mu' \circ (id \otimes \varepsilon) \circ \Delta, \eqno(6\mbox{b})$$
$$m \circ (S \otimes id) \circ \Delta = \varepsilon 
  = m \circ (id \otimes S) \circ \Delta, \eqno(6\mbox{c})$$
where $id$ denotes the identity mapping, 
$\mu : {\cal C} \otimes {\cal A} \longrightarrow {\cal A}, \quad 
  \mu' : {\cal A} \otimes {\cal C} \longrightarrow {\cal A} $ 
are the canonical isomorphisms, defined by 
$\mu(k \otimes a) = ka = \mu'(a \otimes k), \quad \forall a \in {\cal A}, 
  \quad \forall k \in {\cal C}$ 
and $m$ is the multiplication map 
$m : {\cal A} \otimes {\cal A} \longrightarrow {\cal A}, \quad 
  m(a \otimes b) = ab$. 

The multiplication in ${\cal A} \otimes {\cal A}$ follows the rule 
$$(A \otimes B) (C \otimes D) = (-1)^{p(B) p(C)} AC \otimes BC \eqno(7)$$
where $p(X)$ is the $z_2$-grade of $X$, i.e. $p(X) = 0$ for even variables 
and $p(X) = 1$ for odd variables. 

\noindent
{\bf III. LEFT DIFFERENTIAL CALCULUS ON GL$_{p,q}(1\vert 1)$}

In this section, we shall build up the left-invariant differential 
calsulus on the quantum supergroup $GL_{p,q}(1\vert 1)$. The differential 
calculus on the quantum supergroups involves functions on the supergroup, 
differentials and differential forms. It is necessary to point out that, 
to obtain the needed commutation relations for the differential calculus 
we shall not use any specific assumptions. They will be found from natural 
ways. 

\noindent
{\bf A. Left differential algebra}

We first note that the properties of the left exterior differential. 
We can introduce the left exterior differential $\delta_L$ to be a 
${\cal C}$-linear operator that is nilpotent and obeys the graded 
Leibniz rule: 
$$\delta_L^2 = 0, \eqno(8\mbox{a})$$
and 
$$\delta_L(fg) = (fg) \delta_L\hspace*{-0.4cm}^{^{^\leftarrow}} = 
  f(\delta_L g) + (-1)^{p(g)} (\delta_L f) g \eqno(8\mbox{b})$$
where $f$ and $g$ are functions of the group parameters. 

We have seen, in the previous section, that ${\cal A}$ is an associative 
algebra generated by the matrix elements of (1) with the relations (2). 
A differential algebra on ${\cal A}$ is a $z_2$-graded associative algebra 
$\Gamma_L$ equipped with a linear operator $\delta_L $ given (8). Also 
the algebra $\Gamma_L$ has to be generated by 
${\cal A} \cup \delta_L {\cal A}$. 

Firstly, we wish to obtain the relations between the matrix elements of 
$T$ in (1) and their differentials. To do this, we shall use the method of 
ref. 8. So we denote by ${\cal A}_{a\beta}$ the algebra generated by the 
elements $a$ and $\beta$ with the relations 
$$a \beta = q \beta a, \qquad \beta^2 = 0. \eqno(9)$$
If we consider a possible set of commutation relations between generators 
of ${\cal A}_{a\beta}$ and $\delta_L {\cal A}_{\delta_L a \delta_L \beta}$ 
of the form 
$$a ~\delta_L a = A_1 \delta_L a ~a, $$
$$a ~\delta_L \beta = F_{11} \delta_L \beta ~a + F_{12} \delta_L a ~\beta, 
  \eqno(10)$$
$$\beta ~\delta_L a = F_{21} \delta_L a ~\beta + F_{22} \delta_L \beta ~a, $$
$$\beta ~\delta_L \beta = B_1 \delta_L \beta ~\beta, $$
then we can determine the coefficients $A, B$ and $F_{ij}$ in terms of 
the complex deformation parameters $p$, $q$.  To determine them we use 
the consistency of calculus (see, for details, ref. 8). Continuing in 
this way, we can obtain the other relations. The final result is given by 
$$ \delta_L a ~a = pq a~ \delta_L a, $$
$$\delta_L a ~\beta = - q \beta ~\delta_L a + (1 - pq) a~ \delta_L \beta,$$
$$ \delta_La~ \gamma = - p \gamma ~\delta_L a + (1 - pq) a~ \delta_L \gamma, 
  \eqno(11)$$
$$\delta_L a ~d = d ~\delta_L a + (q^{-1} - p) [q p^{-1} 
  \beta~ \delta_L \gamma - \gamma~ \delta_L \beta 
  + (p^{-1} - q) a~ \delta_L d],$$
$$\delta_L \beta~ a = p a~ \delta_L \beta, \qquad 
  \delta_L \beta ~\gamma = p q^{-1} \gamma ~\delta_L \beta + 
  (p - q^{-1}) a~ \delta_L d,$$
$$\delta_L \beta~ \beta = \beta~ \delta_L \beta, \qquad 
  \delta_L \beta~ d = q^{-1} d ~\delta_L \beta + 
  (p^{-1} q^{-1} - 1) \beta~ \delta_L d,$$
$$\delta_L \gamma~ a = q a~ \delta_L \gamma, \qquad 
  \delta_L \gamma ~\beta = q p^{-1} \beta~ \delta_L \gamma + 
   (p^{-1} - q) a~ \delta_L d, $$
$$\delta_L \gamma~ \gamma = \gamma~ \delta_L \gamma, \qquad 
  \delta_L \gamma ~d = p^{-1} d~ \delta_L \gamma + 
   (p^{-1} q^{-1} - 1) \gamma~ \delta_L d, $$
$$\delta_L d ~a = a~ \delta_L d, \qquad 
  \delta_L d ~\beta = - p^{-1} \beta~ \delta_L d, $$
$$\delta_L d~ \gamma = - q^{-1} \gamma~ \delta_L d, \qquad 
  \delta_L d ~d = p^{-1} q^{-1} d~ \delta_L d$$

Note that, for each possible set of the form (10) the emerged equation 
systems admit, of course, at least two solutions. In other words, the 
relations ${\cal A}_{a\beta} - \delta_L {\cal A}_{a\beta}$ appearing 
in eqs. (11) are not unique. For example, the equation 
$$F_{12} F_{22} = 0 = (F_{11} - q A_1) F_{22}, $$
which follows from the consistency of calculus, admits two solutions. 
Here, we choose $F_{22} = 0$. We also note that the coefficients $A_1$ and 
$D_1$ (from $d ~\delta_L d = D_1 \delta_L d ~d$) are, essentialy, 
undetermined. Howover, we have taken them as $A_1 = p^{-1} q^{-1}$ 
and $D_1 = A^{-1}_1$ since these lead to the standard R-matrix 
[see, eq. (69)]. 

To find the commutation relations between differentials, we apply the 
exterior differential $\delta_L$ on the relations (11) and use the nilpotency 
of $\delta_L$. Then it is easy to see that 
$$\delta_L a \delta_L \beta = p^{-1} \delta_L  \beta \delta_L a, 
  \qquad 
  \delta_L d \delta_L \beta = p^{-1} \delta_L \beta \delta_L d, $$
$$\delta_L a \delta_L \gamma = q^{-1} \delta_L \gamma \delta_L a, \qquad 
  \delta_L d \delta_L \gamma = q^{-1} \delta_L \gamma \delta_L d, \eqno(12)$$
$$\delta_L a \delta_L d = - \delta_L d \delta_L a, \qquad 
   (\delta_L a)^2 = 0 = (\delta_L d)^2, $$
$$ \delta_L \beta \delta_L \gamma = p q^{-1} \delta_L \gamma \delta_L \beta + 
   (p - q^{-1}) \delta_L d \delta_L a. $$
These relations are the relations of Gr$_{p,q}(1\vert 1)$ in ref. 12, where 
$\alpha = \delta_L a$, $b = \delta_L \beta$, etc. 

Thus we have constructed the differential algebra 
$\Gamma_L = {\cal A} \cup \delta_L {\cal A}$ of the algebra generated by the 
matrix elements of any matrix in GL$_{p,q}(1\vert 1)$. 

\noindent
{\bf B. Hopf algebra structure on $\Gamma_L$} 

We first note that consistency of a differential calculus with commutation 
relations (1) means that the algebra $\Gamma$ is a graded associative algebra 
generated by the set $\{a, \ldots, d, \delta_L a, \ldots, \delta_L d \}$. 
So, it is sufficient only describe the action of co-maps on the subset 
$\{\delta_L a, \ldots, \delta_L d \}$ which is defined in$^{13}$. 

We consider a map 
$\phi_R : \Gamma_L \longrightarrow \Gamma_L \otimes {\cal A}$ such that 
$$\phi_R \circ \delta_L = (\delta_L \otimes id) \circ \Delta. 
  \eqno(13\mbox{a})$$
Thus, we have 
$$\phi_R(\delta_L a) = \delta_L a \otimes a + \delta_L \beta \otimes \gamma,$$
$$\phi_R(\delta_L \beta) = \delta_L \beta \otimes d + \delta_L a \otimes \beta, 
  \eqno(14\mbox{a})$$
$$\phi_R(\delta_L \gamma) = \delta_L \gamma \otimes a + 
  \delta_L d \otimes \gamma,$$
$$\phi_R(\delta_L d) = \delta_L d \otimes d + \delta_L \gamma \otimes \beta. $$
We now define a map $\Delta_R$ as follows: 
$$\Delta_R(u_1 \delta_L v_1 + \delta_L v_2 u_2) = 
  \Delta(u_1) \phi_R(\delta_L v_1) + \phi_R(\delta_L v_2) \Delta(u_2). 
  \eqno(15)$$
Then it can be checked that the map $\Delta_R$ leaves invariant the relations 
(11) and (12). One can also check that the following identities are satisfied: 
$$(\Delta_R \otimes id) \circ \Delta_R = (id \otimes \Delta) \circ \Delta_R, 
  \qquad (id \otimes \epsilon) \circ \Delta_R = id. \eqno(16\mbox{a})$$
But we do not have a coproduct for the differential algebra because the map 
$\phi_R$ does not gives an analog for the derivation property (8b), 
yet. So we consider another map 
$\phi_L : \Gamma_L \longrightarrow {\cal A} \otimes \Gamma_L$ such that 
$$\phi_L \circ \delta_L = (\tau \otimes \delta_L) \circ \Delta, 
  \eqno(13\mbox{b})$$
where $\tau: \Gamma_L \longrightarrow \Gamma_L$ is the linear map of 
degree zero which gives $\tau(a) = (-1)^{p(a)} a$. 
The action of $\phi_L$ on the generators $\delta_L a$, $\delta_L \beta$, 
$\delta_L \gamma$ and $\delta_L d$ as follows:
$$\phi_L(\delta_L a) = a \otimes \delta_L a - 
  \beta \otimes \delta_L \gamma ,$$
$$\phi_L(\delta_L \beta) = a \otimes \delta_L \beta - 
  \beta \otimes \delta_L d, \eqno(14\mbox{b})$$
$$\phi_L(\delta_L \gamma) = - \gamma \otimes \delta_L a + 
  d \otimes \delta_L \gamma,$$
$$\phi_L(\delta_L d) = d \otimes \delta_L d - \gamma \otimes \delta_L \beta. $$
We define a map $\Delta_L$ with again (15) by replacing $L$ with $R$. 
The map $\Delta_L$ also leaves invariant the relations (11) and (12), 
and the following identities are satisfied: 
$$(id \otimes \Delta_L) \circ \Delta_L = (\Delta \otimes id) \circ \Delta_L, 
  \qquad (\epsilon \otimes id) \circ \Delta_L = id. \eqno(16\mbox{b})$$

Let us define the map $\hat{\Delta}$ as 
$$\hat{\Delta} = \Delta_L + \Delta_R \eqno(17)$$
which will allow us to define the coproduct of the differential algebra. 
We denote the restriction of $\hat{\Delta}$ to the algebra ${\cal A}$ by 
$\Delta$ and the extension of $\Delta$ to the differential algebra $\Gamma_L$ 
by $\hat{\Delta}$. It is possible to interpret the expression 
$$\hat{\Delta}\vert_{\cal A} = \Delta \eqno(18)$$
as the definition of $\hat{\Delta}$ on the matrix elements and (17) as 
the definition of $\hat{\Delta}$ on differentials. 

Note that it is not difficult to verify the following identities: 
$$(\Delta_L \otimes id) \circ \Delta_R = 
  (id \otimes \Delta_R) \circ \Delta_L, \eqno(19\mbox{a})$$
and for all $u \in {\cal A}$ 
$$(\tau \otimes \delta_L) \circ \Delta(u) = \Delta_L(\delta_L u), \qquad 
  (\delta_L \otimes id) \circ \Delta(u) = \Delta_R(\delta_L u). 
  \eqno(19\mbox{b})$$
Note that the coproduct can be interpreted as a (left and right) coaction 
of the quantum supergroup GL$_{p,q}(1\vert 1)$ on the differential forms, 
since the extended algebra $\Gamma_L$ is interpreted as an algebra of 
differential forms on GL$_{p,q}(1\vert 1)$. 

Now let us return Hopf algebra structure of $\Gamma_L$. If we define 
a counit $\hat{\epsilon}$ for the differential algebra as 
$$\hat{\epsilon} \circ \delta_L = \delta_L \circ \epsilon = 0 \eqno(20)$$
and 
$$\hat{\epsilon}\vert_{\cal A} = \epsilon, \qquad 
  \epsilon\vert_\Gamma = \hat{\epsilon} \eqno(21)$$
we have 
$$\hat{\epsilon}(\delta_L a) = \hat{\epsilon}(\delta_L \beta) = 
  \hat{\epsilon}(\delta_L \gamma) = \hat{\epsilon}(\delta_L d) = 0, \eqno(22)$$
where 
$$\hat{\epsilon}(u_1 \delta_L v_1 + \delta_L v_2 u_2) = 
  \epsilon(u_1) \hat{\epsilon}(\delta_L v_1) + 
  \hat{\epsilon}(\delta_L v_2) \epsilon(u_2). \eqno(23)$$
Here we used the fact that $\delta_L (1) = 0$. 

The next step is to obtain a coinverse $\hat{S}$. For this, it suffices to 
define $\hat{S}$ such that 
$$\hat{S} \circ \delta_L = \delta_L \circ S \eqno(24)$$
and 
$$\hat{S}\vert_{\cal A} = S, \qquad 
  S\vert_\Gamma = \hat{S} \eqno(25)$$
where 
$$\hat{S}(u_1 \delta_L v_1 + \delta_L v_2 u_2) = 
  \hat{S}(\delta_L v_1) S(u_1) + S(u_2) \hat{S}(\delta_L v_2). \eqno(26)$$
Thus the action of $\hat{S}$ on the generators $\delta_L a$, $\delta_L \beta$, 
$\delta_L \gamma$ and $\delta_L d$ is as follows: 
$$\hat{S}(\delta_L a) = (- A \delta_L a + B \delta_L \gamma) A - 
  (A \delta_L \beta - B \delta_L d) C, $$
$$\hat{S}(\delta_L \beta) = (- A \delta_L a + B \delta_L \gamma) A - 
  (A \delta_L \beta - B \delta_L d) D, \eqno(27)$$
$$\hat{S}(\delta_L \gamma) = (C \delta_L a - D \delta_L \gamma) A + 
  (C \delta_L \beta - D \delta_L d) C, $$
$$\hat{S}(\delta_L d) = (C \delta_L a - D \delta_L \gamma) B + 
  (C \delta_L \beta - D \delta_L d) D.$$

Note that, it is easy to check that $\hat{\epsilon}$ and $\hat{S}$ leave 
invariant the relations (11) and (12). Consequently, we can say that the 
structure $(\Gamma_L, \hat \Delta, \hat \varepsilon, \hat S)$ is a 
graded Hopf algebra. 

\noindent
{\bf C. The left Cartan-Maurer one-forms in $\Gamma_L$}

As in analogy with the left-invariant one-forms on a Lie group in classical 
differential geometry, one can construct the matrix valued one-form $\Omega_L$ 
where 
$$\Omega_L = \left(\matrix {\theta_1 & u_1 \cr u_2 & \theta_2 \cr }\right) = 
  T^{-1} \delta_L T. \eqno(28)$$
Each element of $\Omega_L$ is left-invariant. For, if $T'$ is any fixed 
element of GL$_{p,q}(1\vert 1)$, the left translation by $T'$ is given by 
$$T \longrightarrow T' T,$$
while
$$(T' T)^{-1} \delta_L (T' T) = T^{-1} \delta_L T.$$
This allows us to make explicit calculations in many important groups. 

If we set 
$$T^{-1} = \left(\matrix{ 
    A & B \cr 
    C & D \cr} \right) \eqno(29)$$
as the superinverse [see, eq. (5)] of $T \in GL_{p,q}(1\vert 1)$, we write 
the matrix elements (left one-forms) of $\Omega_L$ as follows 
$$\theta_1 = A \delta_L a + B \delta_L \gamma, \qquad 
  u_1 = A \delta_L \beta + B \delta_L d, $$
$$\theta_2 = D \delta_L d + C \delta_L \beta, \qquad 
  u_2 = C \delta_L a + D \delta_L \gamma. \eqno(30)$$

In this section, we wish to obtain the commutation relations between the 
generators of ${\cal A}$ and one-forms, and so the relations between 
one-forms. For this reason, we need the commutation relations of the 
matrix elements of $T$ and $T^{-1}$. Some calculations give the commutation 
relations between them as follows: 
$$a A = pq A a + 1 - pq, \qquad d A = A d, $$
$$a D = D a, \qquad d D = pq D d + 1 - pq, $$
$$a B = q B a, \qquad d B = q B d, $$
$$a C = p C a, \qquad d C = p C d, $$
$$\beta A = q A \beta, \qquad \gamma A = p A \gamma, \eqno(31)$$
$$\beta D = q D \beta, \qquad \gamma D = p D \gamma, $$
$$\beta B = B \beta, \qquad \gamma B = - pq B \gamma, $$
$$\beta C = - pq C \beta, \qquad \gamma C = C \gamma.$$

Using these relations, we now find the commutation relations of the matrix 
entries of $T$ with those of $\Omega_L$ : 
$$\theta_1 a = pq a \theta_1 + (pq - 1) \beta u_2, $$
$$\theta_1 \beta = - \beta \theta_1 + (1 - pq) a u_1 - 
  p^{-1} q^{-1} (pq - 1)^2 \beta \theta_2, $$
$$\theta_1 \gamma = - pq \gamma \theta_1 + (1 - pq) d u_2, $$
$$\theta_1 d = d \theta_1 + (pq - 1) \gamma u_1 + 
  p^{-1} q^{-1} (pq - 1)^2 d \theta_2, $$
$$  u_1 a = p a u_1 + (p - q^{-1}) \beta \theta_2,$$
$$u_1 \beta = q^{-1} \beta u_1, \qquad 
  u_1 d = q^{-1} d u_1, $$
$$u_1 \gamma = p \gamma u_1 + (p - q^{-1}) d \theta_2, \eqno(32)$$
$$u_2 a = q a u_2, \qquad 
  u_2 \gamma = q \gamma u_2,$$
$$u_2 \beta = p^{-1} \beta u_2 + (p^{-1} - q) a \theta_2, $$
$$u_2 d = p^{-1} d u_2 + (p^{-1} - q) \gamma \theta_2,$$
$$\theta_2 a = a \theta_2, \qquad 
  \theta_2 \beta = - p^{-1} q^{-1} \beta \theta_2$$
$$\theta_2 \gamma = - \gamma \theta_2, \qquad 
  \theta_2 d = p^{-1} q^{-1} d \theta_2, $$

To obtain commutation relations among the left Cartan-Maurer one-forms, 
we shall use the commutation relations of the matrix elements of $T^{-1}$ 
with the differentials of the matrix elements of $T$ which are given in the 
following 
$$\delta_L a ~A = p^{-1} q^{-1} A \delta_L a + (p^{-1} q^{-1} - 1) 
  [B \delta_L \gamma - C \delta_L \beta + (pq -1) D \delta_L d],$$
$$\delta_L a ~B = - p^{-1} B \delta_L a + (q - p^{-1}) D \delta_L \beta,$$
$$\delta_L a ~C = - q^{-1} C \delta_L a + (p - q^{-1}) D \delta_L \gamma, $$
$$\delta_L a ~D = D \delta_L a, $$
$$\delta_L \beta ~A = p^{-1} A \delta_L \beta + (p^{-1} - q) B \delta_L d, $$
$$\delta_L \beta ~B = q p^{-1} B \delta_L \beta, \qquad 
  \delta_L \beta ~D = q D \delta_L \beta,$$
$$\delta_L \beta ~C = C \delta_L \beta + (1 - p q) D \delta_L d, $$
$$\delta_L \gamma ~A = q^{-1} A \delta_L \gamma + (q^{-1} - p) C \delta_L d,$$
$$\delta_L \gamma ~B = B \delta_L \gamma + (1 - pq) D \delta_L d, $$
$$\delta_L \gamma ~C = p q^{-1} C \delta_L \gamma, \qquad 
  \delta_L \gamma ~D = p D \delta_L \gamma, \eqno(33)$$
$$\delta_L d ~A = A \delta_L d, \qquad 
  \delta_L d ~B = - q B \delta_L d,$$
$$\delta_L d ~C = - p C \delta_L d, \qquad 
  \delta_L d ~D = p q D \delta_L d. $$

Using these relations, we obtain the commutation relations of the left 
Cartan-Maurer forms with the differentials of the matrix elements of $T$ as 
follows:
$$\theta_1 \delta_L a = - \delta_L a ~\theta_1 + 
  (1 - p^{-1} q^{-1}) \delta_L \beta ~u_2, \qquad 
  \theta_1 \delta_L d = - \delta_L d ~\theta_1, $$
$$\theta_1 \delta_L \beta = \delta_L \beta ~\theta_1, \qquad 
  \theta_1 \delta_L \gamma =  \delta_L \gamma ~\theta_1 + 
  (p^{-1} q^{-1} - 1) \delta_L d ~u_2,$$
$$u_1 \delta_L a =  q^{-1} \delta_L a ~u_1 + (p - q^{-1}) 
  \delta_L \beta ~(\theta_1 - \theta_2), \qquad 
  u_1 \delta_L \beta = p \delta_L \beta ~u_1, $$
$$u_1 \delta_L d = p \delta_L d ~u_1, \qquad 
  u_1 \delta_L \gamma =  p \delta_L \gamma ~u_1 + 
  (p - q^{-1}) \delta_L d ~(\theta_1 - \theta_2), \eqno(34)$$
$$u_2 \delta_L a =  p^{-1} \delta_L a ~u_2, \qquad 
  u_2 \delta_L \beta = p^{-1} \delta_L \beta ~u_2, $$
$$u_2 \delta_L \gamma =  p^{-1} \delta_L \gamma ~u_2, \qquad 
  u_2 \delta_L d = p^{-1} \delta_L d ~u_2, $$
$$\theta_2 \delta_L a = - \delta_L a ~\theta_2 + 
  (1 - p^{-1} q^{-1}) \delta_L \beta ~u_2,$$
$$\theta_2 \delta_L \gamma =  \delta_L \gamma ~\theta_2 + 
  (p^{-1} q^{-1} - 1) \delta_L d ~u_2, $$
$$\theta_2 \delta_L \beta =  \delta_L \beta ~\theta_2, \quad 
   \theta_2 \delta_L d = - \delta_L d ~\theta_2. $$

We now obtain the commutation relations of the left Cartan-Maurer forms 
$$u_1 \theta_1 = pq \theta_1 u_1 + (1 - pq) \theta_2 u_1, \qquad 
  u_1 \theta_2 = \theta_2 u_1, $$
$$\theta_1 u_2 = pq u_2 \theta_1 + (1 - pq) u_2 \theta_2, \qquad 
  u_2 \theta_2 = \theta_2 u_2, $$
$$\theta_1^2 = (pq - 1) u_2 u_1, \qquad \theta_2^2 = 0, \eqno(35)$$
$$u_1 u_2 = p q u_2 u_1, \qquad 
  \theta_1 \theta_2 + \theta_2 \theta_1 = (pq - 1) u_2 u_1. $$

Note that one can check that the action of $\delta_L$ on (32) and also 
(34), (35) is consistent. These relations allow to evaluate the Lie algebra of 
GL$_{p,q}(1\vert 1)$ via the generators of Lie algebra related to the left 
one-forms. 

\noindent
{\bf D. Quantum superalgebra}

The commutation relations of Cartan-Maurer forms allow us to construct the 
algebra of the generators. To obtain the quantum Lie superalgebra of the 
Lie algebra generators, using (30) we first write the left Cartan-Maurer 
forms as 
$$\delta_L a = a \theta_1 + \beta u_2, \qquad 
  \delta_L \beta = a u_1 + \beta \theta_2, $$
$$\delta_L d = d \theta_2 + \gamma u_1, \qquad 
  \delta_L \gamma = \gamma \theta_1 + d u_2. 
  \eqno(36)$$
The left differential $\delta_L$ can then the expressed in the form 
$$\delta_L  = T^L_1 \theta_1 + T^L_2 \theta_2 + \nabla^L_+ u_1 + 
  \nabla^L_- u_2.   \eqno(37\mbox{a})$$
Here $T^L_1$, $T^L_2$ and $\nabla^L_\pm$ are the quantum Lie algebra 
generators. We now shall obtain the commutation relations of these generators. 
To do this, let us consider an arbitrary function $f$ of the matrix elements 
of $T$ and write the eq. (37a) as follows: 
$$\delta_L f = (f T^L_i) \theta_i + (f \nabla^L_i) u_i, \eqno(37\mbox{b})$$
where 
$$\theta_i \in \{\theta_1, \theta_2\}, \qquad u_i \in \{u_1, u_2\}, \qquad 
  \nabla^L_i \in \{\nabla^L_-, \nabla^L_+ \}.$$
Using the nilpotency of the left exterior differential $\delta_L$ one has 
$$(f T^L_i) \delta_L \theta_i + (f \nabla^L_i) \delta_L u_i = 
  (f T^L_i) (T^L_j \theta_j + \nabla^L_j u_j) \theta_i - 
  (f \nabla^L_i) (T^L_j \theta_j + \nabla^L_j u_j) u_i. \eqno(38)$$
So we need the four two-forms. To obtain these, using the nilpotency of the 
left differential $\delta_L$, we can write $\delta_L \Omega_L$ of the form 
$$\delta_L \Omega_L = \Omega_L \sigma_3 \Omega_L \sigma_3, \qquad 
  \sigma_3 = \left(\matrix{1 & 0 \cr 0 & - 1 \cr}\right). \eqno(39)$$
In terms of the two-forms, these become 
$$\delta_L \theta_1 = \theta_1^2 - u_1 u_2, \qquad 
  \delta_L u_1 = \theta_1 u_1 - u_1 \theta_2, $$
$$\delta_L \theta_2 = \theta_2^2 - u_2 u_1, \qquad 
  \delta_L u_2 = \theta_2 u_2 - u_2 \theta_1. \eqno(40)$$
We can now write down the Cartan-Maurer equations in our case 
$$\left(\matrix { \delta_L \theta_1 & \delta_L u_1   \cr 
                  \delta_L u_2  & \delta_L \theta_2 \cr}\right) 
  = \left(\matrix{ - u_2 u_1  &  - (\theta_1 - \theta_2) u_1 \cr 
  p^{-1} q^{-1} (\theta_1 - \theta_2) u_2 & - u_2 u_1 \cr}\right). \eqno(41)$$
Using the Cartan-Maurer equations we find the following commutation 
relations for the quantum Lie algebra: 
$$[T^L_1, \nabla^L_+] = - \nabla^L_+ + (1 - pq) T^L_1 \nabla^L_+, $$
$$[T^L_2, \nabla^L_+] = \nabla^L_+ - (1 - pq) T^L_1 \nabla^L_+, $$
$$[T^L_1, \nabla^L_-] = \nabla^L_- - (1 - pq) \nabla^L_- T^L_1, $$
$$[T^L_2, \nabla^L_-] = - \nabla^L_- + (1 - pq) \nabla^L_- T^L_1, \eqno(42)$$
$$\nabla^L_+ \nabla^L_- + p q \nabla^L_- \nabla^L_+ = 
  T^L_1 + T^L_2 + (1 - pq) T^L_1 (T^L_1 - T^L_2), $$
$$[T^L_1, T^L_2] = 0, \qquad (\nabla^L_\pm)^2 = 0 $$
or with new generators 
$X^L = T^L_1 + T^L_2$ and $Y^L = T^L_1 - T^L_2$, 
$$[X^L, \nabla^L_\pm] = 0, \qquad [X^L, Y^L] = 0,\qquad (\nabla^L_\pm)^2 = 0,$$
$$[Y^L, \nabla^L_+] = - 2 \nabla^L_+ + (1 - pq) (X^L + Y^L) \nabla^L_+, 
  \eqno(43)$$
$$[Y^L, \nabla^L_-] = 2 \nabla^L_-  - (1 - pq) \nabla^L_- (X^L + Y^L), $$
$$\nabla^L_+ \nabla^L_- + pq \nabla^L_- \nabla^L_+ = X^L + 
  {{1 - pq}\over 2} (X^L + Y^L) Y^L.$$

We also note that the commutation relations (42) of the Lie algebra 
generators should be consistent with monomials of the matrix elements 
of $T$. To proceed, we must evaluate the commutation relations between 
the generators of Lie algebra and the matrix elements of $T$. The 
commuation relations of the generators with the matrix elements can 
be extracted from the Leibniz rule: 
$$\delta_L (f a) = (f a) \delta_L\hspace*{-0.4cm}^{^{^\leftarrow}} 
= f (\delta_L a) + (\delta_L f) a ~\Longrightarrow~ $$
$$ a (T^L_i \theta_i + \nabla^L_i u_i) = \delta_L a + 
  (T^L_i \theta_i + \nabla^L_i u_i) a,   \eqno(44)$$ 
etc. This yields 
$$a T^L_1 = a + pq T^L_1 a, \qquad 
  a \nabla^L_+ = p \nabla^L_+ a,$$
$$a T^L_2 = T^L_2 a + (p - q^{-1}) \nabla^L_+ \beta,$$
$$a \nabla^L_- = \beta + q \nabla^L_- a + (pq - 1) T^L_1 \beta,$$ 
$$\beta T^L_1 = T^L_1 \beta, \qquad 
  \beta \nabla^L_- = - p^{-1} \nabla^L_- \beta,$$
$$\beta T^L_2 = \beta + p^{-1} q^{-1} T^L_2 \beta + (q - p^{-1}) 
  [\nabla^L_- a + (p - q^{-1}) T^L_1 \beta],$$
$$\beta \nabla^L_+ = a - q^{-1} \nabla^L_+ \beta + (pq - 1) 
  T^L_1 a, \eqno(45)$$ 
$$\gamma T^L_1 = \gamma + pq T^L_1 \gamma, \qquad 
  \gamma \nabla^L_+ = - p \nabla^L_+ \gamma,$$
$$\gamma T^L_2 = T^L_2 \gamma + (q^{-1} - p) \nabla^L_+ d,$$
$$\gamma \nabla^L_- = d - q \nabla^L_- \gamma + (pq - 1) T^L_1 d,$$ 
$$d T^L_1 = T^L_1 d, \qquad 
  d \nabla^L_- = p^{-1} \nabla^L_- d, $$
$$d T^L_2 = d + p^{-1} q^{-1} T^L_2 d + (q - p^{-1}) [(p - q^{-1}) T^L_1 d - 
  \nabla^L_- \gamma],$$ 
$$d \nabla^L_+ = \gamma + q^{-1} \nabla^L_+ d + (pq - 1) T^L_1 \gamma.$$ 

Notice that these commutation relations must be consistent. In fact, 
for example, it is easy to see that the nilpotency of 
$(\nabla^L_\pm)^2$ 
is consistent with 
$$(\nabla^L_+)^2 a = p^{-2} a (\nabla^L_+)^2, \qquad 
  (\nabla^L_-)^2 a = q^{-2} a (\nabla^L_-)^2. $$
Similarly, one can check the other relations. 

\vfill\eject
\noindent
{\bf IV. RIGHT DIFFERENTIAL CALCULUS ON GL$_{p,q}(1\vert 1)$}

In this section, we shall build up the right-invariant differential calculus 
on the quantum supergroup $GL_{p,q}(1\vert 1)$. 

\noindent
{\bf A. Right differential Algebra}

We first note that the properties of the right exterior differential. 
The basic differential operator $\delta_R$ which is linear and satisfies 
the standard properties as follows: 
the nilpotency 
$$ \delta_R^2 = 0 \eqno(46\mbox{a})$$
and the graded Leibniz rule 
$$\delta_R(fg) = (\delta_R f) g + (-1)^{p(f)} f (\delta_R g) 
  \eqno(46\mbox{b})$$
where $f$ and $g$ are functions of the group parameters. 

In Analogy with previous sections let us begin obtaining the 
commutation relations of the group parameters with their differentials. 
Using the method of sec. 3.1 we obtain the commutation relations 
between the generators of ${\cal A}$ and their right differentials 
(the generators of $\delta_R {\cal A}$) as follows: 
$$ a ~\delta_R a = pq \delta_R a~ a, \qquad 
   a~ \delta_R \beta = q \delta_R \beta~ a + (pq - 1) \delta_R a~ \beta, $$
$$ a ~\delta_R \gamma = p \delta_R \gamma~ a + (pq - 1) \delta_R a~ \gamma,$$
$$a~ \delta_R d = \delta_R d~ a + (p - q^{-1}) [\delta_R \gamma~ \beta 
  - q p^{-1} \delta_R \beta~ \gamma + (q - p^{-1}) \delta_R a~ d],$$
$$\beta~ \delta_R a = - p \delta_R a~ \beta, \qquad 
  \beta~ \delta_R \beta = \delta_R \beta~ \beta, $$
$$ \beta~ \delta_R \gamma = pq^{-1} \delta_R \gamma~ \beta + 
   (p - q^{-1}) \delta_R a~ d, $$
$$ \beta ~\delta_R d = - q^{-1} \delta_R d~ \beta + 
   (1 - p^{-1} q^{-1}) \delta_R \beta~ d,$$
$$\gamma~ \delta_R a = - q \delta_R a~ \gamma, \qquad 
  \gamma ~\delta_R \gamma = \delta_R \gamma~ \gamma,$$
$$\gamma~ \delta_R \beta = q p^{-1} \delta_R \beta~ \gamma + 
   (p^{-1} - q) \delta_R a~ d,   \eqno(47)$$
$$\gamma~ \delta_R d = - p^{-1} \delta_R d~ \gamma + 
  (1 - p^{-1} q^{-1}) \delta_R \gamma~ d,$$
$$ d~ \delta_R a = \delta_R a~ d, \qquad 
   d~ \delta_R \beta = p^{-1} \delta_R \beta~ d,$$
$$ d~ \delta_R \gamma = q^{-1} \delta_R \gamma~ d, \qquad 
   d~ \delta_R d = p^{-1} q^{-1} \delta_R d~ d.$$ 

To find the commutation relations between differentials, we apply the 
exterior differential $\delta_R$ on the relations (47) and use the nilpotency 
of $\delta_R$. They are the same with (12), as expected. 

Thus we have constructed the right differential algebra 
$\Gamma_R = {\cal A} \cup \delta_R {\cal A}$ of the algebra generated by the 
matrix elements of any matrix in GL$_{p,q}(1\vert 1)$. Again,the 
structure $(\Gamma_R, \hat \Delta, \hat \varepsilon, \hat S)$ is a 
graded Hopf algebra provided that subscript $R$ replacing with $L$ in 
sec. 3.B.

\noindent
{\bf B. The right Cartan-Maurer one-forms in $\Gamma_R$}

As in analogy with the right-invariant one-forms on a Lie group in classical 
differential geometry, one can construct the matrix valued one-form $\Omega_R$ 
where 
$$\Omega_R = \delta_R T ~T^{-1}. \eqno(48)$$
Then we write the matrix elements (right one-forms) of $\Omega_R$ as follows 
$$w_1 = \delta_R a A + \delta_R \beta C, \qquad 
  v_1 = \delta_R a B + \delta_R \beta D, $$
$$w_2 = \delta_R \gamma B + \delta_R d D, \qquad 
  v_2 = \delta_R \gamma A + \delta_R d C. \eqno(49)$$

The commutation relations of the matrix elements of $T$ and $T^{-1}$ are 
given by (31). Using (31), we now find the commutation relations of the 
matrix entries of $T$ with those of $\Omega_R$ : 
$$a w_1 = pq w_1 a, \qquad 
  a v_1 = q v_1 a ,$$
$$a v_2 = p v_2 a + (p - q^{-1}) w_1 \gamma, $$
$$a w_2 = w_2 a + p q^{-1} (q - p^{-1})^2 w_1 a + 
  (p^{-1} q^{-1} - 1) v_1 \gamma, $$
$$\beta w_1 = - pq w_1 \beta, \qquad 
  \beta v_1 = q v_1 \beta, $$
$$\beta v_2 = p v_2 \beta + (p - q^{-1}) w_1 d, \eqno(50)$$
$$\beta w_2 = - w_2 \beta - p q^{-1} (q - p^{-1})^2 w_1 \beta + 
  (1 - p^{-1} q^{-1}) v_1 d, $$
$$\gamma w_1 = - w_1 \gamma, \qquad 
  \gamma v_1 = p^{-1} v_1 \gamma + (p^{-1} - q) w_1 a,$$
$$\gamma v_2 = q^{-1} v_2 \gamma, \qquad 
  \gamma w_2 = - p^{-1} q^{-1} w_2 \gamma + (1 - p^{-1} q^{-1}) v_2 a,$$
$$d w_1 = w_1 d, \qquad d v_2 = q^{-1} v_2 d, $$
$$  d v_1 = p^{-1} v_1 d + (p^{-1} - q) w_1 \beta, $$
$$d w_2 = p^{-1} q^{-1} w_2 d + (p^{-1} q^{-1} - 1) v_2 \beta.$$

To obtain the commutation relations among the right Cartan-Maurer 
one-forms, we use the commutation relations of the matrix elements 
of $T^{-1}$ with the differentials of the matrix elements of $T$ 
which are given in the following 
$$A \delta_R a = p^{-1} q^{-1} \delta_R a A, \qquad 
  A \delta_R d = \delta_R d A, $$
$$A \delta_R \beta = q^{-1} \delta_R \beta A, \qquad 
  A \delta_R \gamma = p^{-1} \delta_R \gamma A, $$
$$D \delta_R a = \delta_R a D, $$
$$D \delta_R d = pq \delta_R d D + (1 - pq) [\delta_R \beta C - 
  \delta_R \gamma B + (1 - p^{-1} q^{-1}) \delta_R a A], $$
$$D \delta_R \beta = p \delta_R \beta D + (p - q^{-1}) \delta_R a B, $$
$$D \delta_R \gamma = q \delta_R \gamma D + (p - p^{-1}) \delta_R a C, $$
$$B \delta_R a = - q^{-1} \delta_R a B, \qquad 
  B \delta_R \gamma = \delta_R \gamma B + (p^{-1} q^{-1} - 1) \delta_R a A, 
  \eqno(51)$$
$$B \delta_R \beta = pq^{-1} \delta_R \beta B, \qquad 
  B \delta_R d = - p \delta_R d B + (q^{-1} - p) \delta_R \beta A$$
$$C \delta_R a = - p^{-1} \delta_R a C, \qquad 
  C \delta_R \beta = \delta_R \beta C + (1 - p^{-1} q^{-1}) \delta_R a A, $$
$$C \delta_R \gamma = qp^{-1} \delta_R \gamma C, \qquad 
  C \delta_R d = - q \delta_R d C + (p^{-1} - q) \delta_R \gamma A.$$

Using these relations, we obtain the commutation relations of the right 
Cartan-Maurer forms with the differentials of the matrix elements of $T$ as 
follows:
$$w_1 \delta_R a = - \delta_R a ~w_1, \qquad 
  w_1 \delta_R \gamma = \delta_R \gamma ~w_1,$$
$$w_1 \delta_R \beta = \delta_R \beta ~w_1 + 
  (1 - p^{-1} q^{-1}) \delta_R a ~v_1,$$
$$w_1 \delta_R d = - \delta_R d ~w_1 + 
  (p^{-1} q^{-1} - 1) \delta_R \gamma ~v_1,$$
$$v_1 \delta_R a =  q^{-1} \delta_R a ~v_1, \qquad 
  v_1 \delta_R \beta = q^{-1} \delta_R \beta ~v_1, $$
$$v_1 \delta_R \gamma =  q^{-1} \delta_R \gamma ~v_1, \qquad 
  v_1 \delta_R d = q^{-1} \delta_R d ~v_1, \eqno(52)$$
$$v_2 \delta_R a = q \delta_R a ~v_2, \qquad 
  v_2 \delta_R \gamma =  q \delta_R \gamma ~v_2, $$
$$v_2 \delta_R \beta = q \delta_R \beta ~v_2 + 
  (q - p^{-1}) \delta_R a ~(w_2 - w_1),$$
$$v_2 \delta_R d = q \delta_R d ~v_2 + (q - p^{-1}) \delta_R \gamma ~w_1,$$
$$w_2 \delta_R a = - \delta_R a ~w_2, \qquad 
  w_2 \delta_R \beta = \delta_R \beta ~w_2 + 
  (1 - p^{-1} q^{-1}) \delta_R a ~v_1,$$
$$w_2 \delta_R \gamma =  \delta_R \gamma ~w_2 + 
  p^{-1} q^{-1} (pq - 1)^2 \delta_R \gamma ~w_1, $$
$$w_2 \delta_R d = - pq \delta_R d ~w_2 + (p q - 1) \delta_R d ~w_1 + 
  p^{-1} q^{-1} (pq - 1)^2 \delta_R \gamma ~v_1. $$

We now obtain the commutation relations of the right Cartan-Maurer forms 
$$w_1 v_1 = v_1 w_1, \qquad v_1 w_2 = pq w_2 v_1 + (1 - pq) v_1 w_1, $$
$$w_1 v_2 = v_2 w_1, \qquad w_2 v_2 = pq v_2 w_2 + (1 - pq) w_1 v_2, $$
$$w_1^2 = 0, \qquad w_2^2 = (1 - pq) v_2 v_1, \eqno(53)$$
$$v_1 v_2 = p q v_2 v_1, \qquad w_1 w_2 + w_2 w_1 = (1 - pq) v_2 v_1. $$

Note that one can check that the action of $\delta_R$ on (50), (52) and also 
(53) is consistent. These relations allow us to evaluate the Lie algebra of 
GL$_{p,q}(1\vert 1)$ by relating the generators of the Lie algebra to the 
right one-forms. 

\noindent
{\bf C. Quantum superalgebra}

The commutation relations of Cartan-Maurer forms allow us to construct the 
algebra of the generators. To obtain the quantum Lie superalgebra of the 
Lie algebra generators we first write the Cartan-Maurer forms as 
$$\delta_R a = w_1 a + v_1 \gamma, \qquad 
  \delta_R \beta = w_1 \beta + v_1 d, $$
$$\delta_R d = w_2 d + v_2 \beta, \qquad 
  \delta_R \gamma = w_2 \gamma + v_2 a.   \eqno(54)$$
The differential $\delta_R$ can then the expressed in the form 
$$\delta_R  = w_1 T_1 + w_2 T_2 + v_1 \nabla_+ + v_2 \nabla_-. \eqno(55)$$
Here $T_1$, $T_2$ and $\nabla_{\pm}$ are the quantum Lie algebra generators. 
We now shall obtain the commutation relations of these generators. 
Considering an arbitrary function $f$ of the matrix elements of $T$ and using 
the nilpotency of the exterior differential $\delta_R$ one has 
$$(\delta_R w_i) T_i f + (\delta_R v_i) \nabla_i f = 
  w_i \delta_R T_i f - v_i \delta_R \nabla_i f, \eqno(56)$$
where 
$$w_i \in \{w_1,w_2\}, \qquad v_i \in \{v_1,v_2\}, \qquad 
  \nabla_i \in \{\nabla_+, \nabla_-\}.$$
So we need the four two-forms. To obtain these, using the nilpotency of the 
differential $\delta_R$, we can write $\delta_R \Omega_R$ of the form 
$$\delta_R \Omega_R = \sigma_3 \Omega_R \sigma_3 \Omega_R, \qquad 
  \sigma_3 = \left(\matrix{1 & 0 \cr 0 & - 1 \cr}\right) \eqno(57)$$
In terms of the two-forms, these become 
$$\delta_R w_1 = w_1^2 - v_1 v_2, \qquad 
  \delta_R v_1 = w_1 v_1 - v_1 w_2, $$
$$\delta_R w_2 = w_2^2 - v_2 v_1, \qquad 
  \delta_R v_2 = w_2 v_2 - v_2 w_1. \eqno(58)$$
We can now write down the Cartan-Maurer equations in our case 
$$\left(\matrix { \delta_R w_1 & \delta_R v_1   \cr 
                  \delta_R v_2   & \delta_R w_2 \cr}\right) 
  = \left(\matrix{ - v_1 v_2  & p q (w_1 - w_2) v_1 \cr 
                   - (w_1 - w_2) v_2 & - v_1 v_2 \cr}\right). \eqno(59)$$
Using the Cartan-Maurer equations we find the following commutation 
relations for the quantum Lie algebra: 
$$[T_1,\nabla_+] = - pq \nabla_+ + (p q - 1) T_2 \nabla_+, $$
$$[T_2,\nabla_+] = pq \nabla_+ - (pq - 1) T_2 \nabla_+, $$
$$[T_1,\nabla_-] = pq \nabla_- - (pq - 1) \nabla_- T_2, $$
$$[T_2,\nabla_-] = - pq\nabla_- + (pq - 1) \nabla_- T_2, \eqno(60)$$
$$[T_1,T_2] = 0, \qquad \nabla_{\pm}^2 = 0$$
$$\nabla_- \nabla_+ + p^{-1} q^{-1} \nabla_+ \nabla_- = T_1 + T_2 + 
  (p^{-1} q^{-1} - 1) (T_2^2 + T_1 T_2) $$
or with new generators $X = T_1 + T_2$ and $Y = T_1 - T_2$, 
$$[X,\nabla_{\pm}] = 0, \qquad [X,Y] = 0, \qquad \nabla_{\pm}^2 = 0, $$
$$[Y,\nabla_+] = - 2 pq \nabla_+ + (pq - 1) (X - Y)\nabla_+, $$
$$[Y,\nabla_-] = 2 pq \nabla_- - (pq - 1) \nabla_- (X - Y), \eqno(61)$$
$$\nabla_+ \nabla_- + pq \nabla_- \nabla_+ = 
   pq X + {{1 - pq}\over 2} (X^2 - XY). $$

The commutation relations (60) of the Lie algebra generators should 
be consistent with monomials of the matrix elements of $T$. 
To do this, we evaluate the commutation relations between the 
generators of Lie algebra and the matrix elements of $T$. The commuation 
relations of the generators with the matrix elements can be extracted from 
the Leibniz rule: 
$$\delta_R (a f) = (\delta_R a) f + a (\delta_R  f) \Longrightarrow 
  (w_i T_i + v_i \nabla_i) a = \delta_R a + a (w_i T_i + v_i \nabla_i), 
  \eqno(62)$$
etc. This yields 
$$T_1 a  = a + pq aT_1 + (p - q^{-1}) [(q - p^{-1}) a T_2 + 
  \gamma \nabla_-], $$
$$T_1 \beta  = \beta + pq \beta T_1 + (p - q^{-1}) [(q - p^{-1}) \beta T_2 - 
  d \nabla_-], $$
$$T_1 \gamma  = \gamma T_1 + (p^{-1} - q) a \nabla_+,$$
$$T_1 d  = d T_1 + (p^{-1} - q) \beta \nabla_+, $$
$$T_2 a  = a T_2, \qquad T_2 \gamma  = \gamma + p^{-1} q^{-1} \gamma T_2,$$
$$T_2 \beta  = \beta T_2,  \qquad 
    T_2 d = d + p^{-1} q^{-1} d T_2, \eqno(63)$$
$$\nabla_+ a = \gamma + q a \nabla_+ + (p^{-1} q^{-1} - 1) \gamma T_2, $$
$$\nabla_+ \beta = d - q \beta \nabla_+ + (p^{-1} q^{-1} - 1) d T_2,$$
$$\nabla_+ \gamma = - p^{-1} \gamma \nabla_+, \qquad 
  \nabla_+ d = p^{-1} d \nabla_+, $$
$$\nabla_- a = p a \nabla_-, \qquad 
  \nabla_- \beta = - p \beta \nabla_-, $$
$$\nabla_- \gamma = a -  q^{-1} \gamma \nabla_- + (p^{-1} q^{-1} - 1) a T_2, $$
$$\nabla_- d = \beta + q^{-1} d \nabla_- + (p^{-1} q^{-1} - 1) \beta T_2. $$

Notice that these commutation relations must be consistent. In fact, for 
example, it is easy to see that the nilpotency of $\nabla_\pm^2$ is 
consistent with 
$$\nabla_-^2 a = p^2 a \nabla_-^2, \qquad 
  \nabla_+^2 a = q^2 a \nabla_+^2. $$
Similarly, one can check the other relations. 

\noindent
{\bf V. R-MATRIX APPROACH}

In this section we wish to obtain the relations (11), (12), (32), (34) and 
(35) with the help of a matrix $\hat{R}$ that acts on the square tensor space 
of the supergroup. Of course, the matrix $\hat R$ is a solution of the 
quantum (graded) braided group equation. 

We first consider the quantum superplane and its dual$^3$. The quantum 
superplane $A_p$ is generated by coordinates $x$ and $\theta$, and the 
commutation rules 
$$x \theta = p \theta x, \qquad \theta^2 = 0. \eqno(64)$$
The quantum (dual) superplane $A^*_q$ is generated by coordinates $\varphi$ 
and $y$, and the commutation rules 
$$\varphi^2 = 0, \qquad \varphi y = q^{-1} y \varphi. \eqno(65)$$
We demand that relations (64), (65) are preserved under the action of $T$, 
as a linear transformation, on the quantum superplane and its dual: 
$$T : A_p \longrightarrow A_p, \qquad T : A^*_q \longrightarrow A^*_q. 
  \eqno(66)$$
Let $X = (x, \theta)^t$ and $\hat{X} = (\varphi, y)^t$. Then, as a 
consequence of (66) the points $TX$ and $T\hat{X}$ should belong to $A_p$ 
and $A_q^*$, respectively, which give the relations (2). 

Similarly, let us consider linear transformations $\delta_L T$ with the 
following properties 
$$\delta_L T : A_p \longrightarrow A^*_q, \qquad 
  \delta_L T : A^*_q \longrightarrow A_p.   \eqno(67)$$
Then the points $(\delta_L T) X$ and $(\delta_L T) \hat{X}$ should belong to 
$A_q^*$ and $A_p$, respectively. This case is equivalent to (12). 

Note that the relations (64) can be written as follows 
$$X \otimes X = q^{-1} \hat{R} X \otimes X, \eqno(68)$$
where 
$$ \hat{R} = \left(\matrix{ 
   q &    0       & 0 & 0 \cr
   0 & q - p^{-1} & 1 & 0 \cr
   0 &  qp^{-1}   & 0 & 0 \cr
   0 &    0       & 0 & - p^{-1} \cr} \right). \eqno(69)$$
We can also write mixed relations between the component of $X$ and $\hat X$ 
as follows: 
$$(-1)^{p(X)} X \otimes \hat{X} = p \hat{R} \hat{X} \otimes X, \eqno(70)$$
where $\hat{X} = \delta_L X$. 

Using (66) together with (68) and (70), we now derive anew the quantum 
supergroup relations (2) from the equation 
$$\hat {R} T_1 T_2 = T_1 T_2 \hat {R}, \eqno(71)$$
where, in usual grading tensor notation, $T_1 = T \otimes I$ and 
$T_2 = I \otimes T$. Similarly using (70), we obtain the following equation 
$$ T_1 \delta_L (T_2) = q p^{-1} \hat R^{-1} \delta_L T_1 T_2' \hat R^{-1}, 
  \qquad 
   T' = (-1)^{p(T)} T  \eqno(72\mbox{a})$$
which is equivalent to the relations (11). 
Note that 
$$\delta_L (T_2) = - (\delta_L T)_2. \eqno(73)$$
So, the equation (72a) can be written as 
$$ T_1 (\delta_L T)_2 = - q p^{-1} \hat R^{-1} \delta_L T_1 T_2' \hat R^{-1}.
  \eqno(72\mbox{b})$$
Applying the left exterior differential $\delta_L$ on boht side of (72) one 
has  
$$\delta_L T_1 [\delta_L (T_2)]' = 
  q p^{-1} \hat R^{-1} \delta_L T_1 \delta_L (T'_2) \hat R^{-1}, 
  \eqno(74)$$
which gives the relations (12). Taking $\delta_L T = T \Omega_L$ and using 
(72) one obtains 
$$\Omega_1 T'_2 = - p q^{-1} T_2 \hat R \Omega_2 \hat R, \eqno(75)$$
which gives the relations (34). Finally, from (74) we find that 
$$\Omega_1 [\delta_L (T_2)]' = \delta_L (T_2) \hat R \Omega'_2 \hat R^{-1} 
 \eqno(76)$$
and from (76) 
$$\hat R \Omega_2 \hat R \Omega'_2 = - qp^{-1} 
  \Omega_2 \hat R \Omega'_2 \hat R^{-1}. \eqno(77)$$
These equations are equivalent to (34) and (35), respectivelly. Similar 
formulas can be also derived for the right commutation relations. 

\noindent
{\bf VI. CONCLUSION}

To conclude, we introduce here commutation relations between the group 
parameters and their partial derivatives and thus illustrate the connection 
between the relations in Sec. IIID and in Sec. IVC, and the relations which 
will be now obtained. 

To proceed, let us first obtain the relations of the group parameters 
with their partial derivatives. We know that the right exterior 
differential $\delta_R$ can be expressed of the form 
$$\delta_R f = 
  (\delta_R a \partial_a + \delta_R \beta \partial_\beta + 
  \delta_R \gamma \partial_\gamma + \delta_R d \partial_d)f. \eqno(78)$$
Then, replacing $f$ with $af$, etc. we obtain the following commutation 
relations 
$$\partial_a a = 1 + pq a \partial_a + (pq - 1) 
  [(1 - p^{-1} q^{-1}) d \partial_d + \beta \partial_\beta + 
  \gamma \partial_\gamma], $$
$$\partial_a \beta = p \beta \partial_a + (q^{-1} - p) d \partial_\gamma, $$
$$\partial_a \gamma = q \gamma \partial_a + (q - p^{-1}) d \partial_\beta, 
  \qquad \partial_a d = d \partial_a, $$
$$\partial_\beta a = q a \partial_\beta + (p^{-1} - q) \gamma \partial_d, 
   \qquad \partial_\beta d = p^{-1} d \partial_\beta, $$
$$\partial_\beta \beta = 1 - \beta \partial_\beta + 
  (p^{-1} q^{-1} - 1) d \partial_d, \qquad 
  \partial_\beta \gamma = - q p^{-1} \gamma \partial_\beta, $$
$$\partial_\gamma a = p a \partial_\gamma + (p - q^{-1}) \beta \partial_d, 
  \qquad \partial_\gamma \beta = - p q^{-1} \beta \partial_\gamma, $$
$$\partial_\gamma \gamma = 1 - \gamma \partial_\gamma + 
  (p^{-1} q^{-1} - 1) d \partial_d, \qquad 
  \partial_\gamma d = q^{-1} d \partial_\gamma, \eqno(79)$$
$$\partial_d a = a \partial_d,  \qquad 
  \partial_d \beta = q^{-1} \beta \partial_d, $$
$$\partial_d \gamma = p^{-1} \gamma \partial_d, \qquad 
  \partial_d d = 1 + p^{-1} q^{-1} d \partial_d.$$
We thus find the commutation relations between the derivatives. 
These relations can be obtained by using the nilpotency of the right exterior 
differential $\delta_R$ and they have the form 
$$\partial_a \partial_\beta = p^{-1} \partial_\beta \partial_a, \qquad 
  \partial_d \partial_\beta = p^{-1} \partial_\beta \partial_d, $$
$$\partial_a \partial_\gamma = q^{-1} \partial_\gamma \partial_a, \qquad 
  \partial_d \partial_\gamma = q^{-1} \partial_\gamma \partial_d, $$
$$\partial_\beta \partial_\gamma = - p q^{-1} \partial_\gamma \partial_\beta, 
  \qquad  \partial_\beta^2 = 0 =    \partial_\gamma^2, \eqno(80)$$
$$ \partial_a \partial_d = \partial_d \partial_a + 
  (p - q^{-1}) \partial_\gamma \partial_\beta. $$
The (graded) Hopf algebra structure for $\partial$ is given by 
$$\Delta(\partial_a) = \partial_a \otimes \partial_a + 
  \partial_\beta \otimes \partial_\gamma, \qquad 
  \Delta(\partial_\beta) = \partial_a \otimes \partial_\beta + 
  \partial_\beta \otimes \partial_d, $$
$$\Delta(\partial_d) = \partial_d \otimes \partial_d + 
  \partial_\gamma \otimes \partial_\beta, \qquad 
  \Delta(\partial_\gamma) = \partial_\gamma \otimes \partial_a + 
  \partial_d \otimes \partial_\gamma, \eqno(81)$$
$$\varepsilon(\partial_a) = 1 = \varepsilon(\partial_d), \qquad 
  \varepsilon(\partial_\beta) = 0 = \varepsilon(\partial_\gamma), $$
$$S(\partial_a) = \partial_a^{-1} + \partial_a^{-1}  \partial_\beta 
   \partial_d^{-1} \partial_\gamma \partial_a^{-1}, \qquad 
  S(\partial_\beta) = - \partial_a^{-1} \partial_\beta \partial_d^{-1}, $$
$$S(\partial_d) = \partial_d^{-1} + \partial_d^{-1}  \partial_\gamma 
   \partial_a^{-1} \partial_\beta \partial_d^{-1}, \qquad 
  S(\partial_\gamma) = - \partial_d^{-1} \partial_\gamma \partial_a^{-1}, $$
provided that the formal inverses $\partial_a^{-1}$ and $\partial_d^{-1}$ 
exist. However these co-maps do not leave invariant the relations (79). 

We know, from Sec. IVC, that the right exterior differential $\delta_R$ can 
be expressed in the form (55), which we repeat here, 
$$\delta_R f = (w_1 T_1 + v_1 \nabla_+ + v_2 \nabla_- + w_2 T_2) f. \eqno(82)$$
Considering (78) together (82) and using (54) one has 
$$T_1 = a \partial_a + \beta \partial_\beta, \qquad 
  \nabla_+ = \gamma \partial_a + d \partial_\beta, $$
$$T_2 = d \partial_d + \gamma \partial_\gamma, \qquad 
  \nabla_- = a \partial_\gamma + \beta \partial_d. \eqno(83)$$
Using the relations (79), (80) one can check that the relations of the 
generators in (83) coincide with (60). It can also be verified that, 
the action of the generators in (83) on the group parameters coincide 
with (63). 

Above, we noted that the right exterior differential must be placed 
before the right derivatives. With a similar consideration, we can say 
that the left exterior differential must be placed after the left derivatives: 
$$\delta_L f = f (\partial_a^L\delta_R a + \partial_\beta^L \delta_R \beta 
  + \partial_\gamma^L \delta_R \gamma + \partial_d^L \delta_R d ). \eqno(84)$$
Replacing $f$ with $f a$, etc. we get the following relations 
$$a \partial_a^L a = 1 + pq \partial_a^L a, \qquad 
  a \partial_d^L = \partial_d^L a, $$
$$a \partial_\beta^L = p \partial_\beta^L a, \qquad 
  a \partial_\gamma^L = q \partial_\gamma^L a, $$
$$\beta \partial_a^L = q \partial_a^L \beta, \qquad 
  \beta \partial_\beta^L = 1 - \partial_\beta^L \beta + 
  (pq - 1) \partial_a^L a, $$
$$\beta \partial_\gamma^L = - q p^{-1} \partial_\gamma^L \beta, \qquad 
  \beta \partial_d^L = p^{-1} \partial_d^L \beta + 
  (q - p^{-1}) \partial_\gamma^L a, $$
$$\gamma \partial_a^L = p \partial_a^L \gamma, \qquad 
  \gamma \partial_\gamma^L = 1 - \partial_\gamma^L \gamma + 
  (pq - 1) \partial_a^L a, $$
$$\gamma \partial_\beta^L = - p q^{-1} \partial_\beta^L \gamma, \qquad 
 \gamma \partial_d^L = q^{-1} \partial_d^L \gamma + 
  (q^{-1} - p) \partial_\beta^L a,$$
$$d \partial_a^L = \partial_a^L d, \qquad 
  d \partial_\beta^L = q^{-1} \partial_\beta^L d + 
  (p - q^{-1}) \partial_a^L \gamma, \eqno(85)$$
$$d \partial_\gamma^L = p^{-1} \partial_\gamma^L d + 
  (p^{-1} - q) \partial_a^L \beta,$$
$$d \partial_d^L = 1 + p^{-1} q^{-1} \partial_d^L d + (1 - p^{-1} q^{-1}) 
  [(1 - pq) \partial_a^L a + \partial_\beta^L \beta + 
  \partial_\gamma^L \gamma].$$
The commutation relations between the left derivatives are the same with (80). 

Finally, expressing the left exterior differential of the form (37b) and 
comparing (84) by help of (36) we have 
$$T_1^L = \partial_a^L a+ \partial_\gamma^L \gamma, \qquad 
  \nabla_+^L = \partial_\beta^L a + \partial_d^L \gamma, $$
$$T_2^L = \partial_\beta^L \beta + \partial_d^L d, \qquad 
  \nabla^L_- = \partial_a^L \beta + \partial_\gamma^L d. \eqno(86)$$
Using these with (85) one can check the relations (42) and (45). 

\noindent
{\bf VII. DISCUSSION} 

The starting point of the present paper is to evaluate the 
$(p,q)$-commutation relations of the matrix elements with their 
differentials. Later, using these relations the 
$(p, q)$-commutation relations of the matrix elements with 
the Cartan-Maurer forms are obtained without any further 
assumptions. The commutation relations of the Cartan-Maurer 
forms are not obtained by using the $(p,q)$-commutation 
relations of the matrix elements with the Cartan-Maurer forms, 
i.e., to obtain the desired commutation relations we have not 
applied the exterior differential $\delta$ on the relations 
of the matrix elements with the Cartan-Maurer forms. Applying 
the exterior differential $\delta$ on the 
relations of the matrix elements to the Cartan-Maurer forms, 
gives the required objects. In this work we have derived the 
$(p,q)$-commutation relations between the matrix elements and 
their differentials without considering an R-matrix at first. 
However we later showed that these relations can also be derived 
using an R-matrix. 

\noindent
{\bf ACKNOWLEDGEMENT}

This work was supported in part by T. B. T. A. K. the 
Turkish Scientific and Technical Research Council. 

\noindent
$^1$ N. Y. Reshetikhin, L. A. Takhtajan, and L. D. Faddeev, Leningrad Math. J. {\bf 1}, 193 (1990); \\
E. Corrigan, D. Fairlie, P. Fletcher, and R. Sasaki, J. Math. Phys. {\bf 31}, 776 (1990). \\
$^2$ L. Alvarez-Gaume, C. Gomes, and G. Sierra, Nucl. Phys. B {\bf 319}, 155 (1989); \\
T. Curtright, D. Fairlie, and C. Zachos, "Quantum groups", in Proc. Argonne Workshop (World Scientific) (1990); \\
D. Fairlie, and C. Zachos, Phys. Lett. B {\bf 256}, 43 (1991). \\
$^3$ Yu I. Manin, Commun. Math. Phys. {\bf 123}, 163 (1989). \\
$^4$ S. L. Woronowicz, Commun. Math. Phys. {\bf 122}, 125 (1989). \\
$^5$ S. L. Woronowicz, Kyoto Univ. {\bf 23}, 117 (1987). \\
$^6$ A. Schirmacher, J. Wess and B. Zumino, Z. Phys. C {\bf 49}, 317 (1990).\\
$^7$ P. Aschieri, and L. Castellani, Int. J. Mod. Phys. A {\bf 8}, 1667 (1993); \\
B. Jurco, Lett. Math. Phys. {\bf 22}, 177 (1991);\\
A. Sudbery, Phys. Lett. B {\bf 284}, 61 (1992);\\
F. M\"uller-Hoissen, J. Phys. A: Math. Gen. {\bf 25}, 1703 (1992). \\
$^8$ S. \c Celik, and S. A. \c Celik, J. Phys. A: Math. Gen. {\bf 31}, 9685 (1998). \\
$^9$ W. Schmidke, S. Vokos, and B. Zumino, Z. Phys. C {\bf 48}, 249 (1990). \\
$^{10}$ S. \c Celik, and S. A. \c Celik, Balkan Phys. Lett. {\bf 3}, 188 (1995).\\
$^{11}$ L. Dabrowski, and L. Wang, Phys. Lett. B {\bf 266}, 51 (1991). \\
$^{12}$ S. \c Celik, J. Math. Phys. {\bf 37}, 3568 (1996). \\
$^{13}$ S. \c Celik, J. Math. Phys. {\bf 40}, 2494 (1999).

\end{document}